\documentclass[11pt]{amsart}

\usepackage{latexsym}
\usepackage{amssymb}
\usepackage{amsmath}
\usepackage{amsfonts}
\usepackage[mathscr]{eucal}
\usepackage{auto-pst-pdf}
\usepackage{times}

\usepackage{amscd}

\usepackage{graphicx}
\usepackage{epsfig}

\usepackage{amsthm}
\usepackage[all,cmtip]{xy}
\usepackage{hyperref}
\usepackage{pst-grad} 
\usepackage{pst-plot} 
\usepackage{comment}
\usepackage{psfrag}
\usepackage{verbatim}

\newtheorem{thm}{Theorem}[section]

\theoremstyle{definition}

\theoremstyle{remark}
\newtheorem*{ack}{Acknowledgements}

\numberwithin{equation}{section}

\DeclareMathOperator{\symrank}{symrank}
\DeclareMathOperator{\Iso}{Isom}

\newcommand{\RP}{\mathbb{R}P}
\newcommand{\CP}{\mathbb{C}P}

\newcommand{\sphere}{\mathrm{\mathbb{S}}}

\newcommand{\SO}{\mathrm{SO}}
\newcommand{\SU}{\mathrm{SU}}

\begin{document}



\title[Symmetry rank and quasipositive curvature]{A note on maximal symmetry rank, quasipositive curvature, and low dimensional manifolds}

\author[F.\ Galaz-Garcia]{Fernando Galaz-Garcia$^\ast$}

\thanks{$^\ast$ The author is part of SFB 878: \emph{Groups, Geometry \& Actions}, at the University of M\"unster.}

\address{Mathematisches Institut, Westf\"alische Wilhelms-Universit\"at M\"unster, Germany}
\email{f.galaz-garcia@uni-muenster.de}

\date{\today}




\begin{abstract}
We show that any effective isometric torus action of maximal rank on a compact Riemannian manifold with positive (sectional) curvature and maximal symmetry rank, that is, on a positively curved sphere, lens space, complex or real projective space, is equivariantaly diffeomorphic to a linear action. We show that a compact, simply connected  Riemannian $4$- or $5$-manifold of quasipositive curvature and maximal symmetry rank must be diffeomorphic to the $4$-sphere, complex projective plane or the $5$-sphere.    \end{abstract}

\maketitle



%
%

\section{Introduction and main results}

The topological classification of compact Riemannian $n$-manifolds  with positive or nonnegative (sectional) curvature is a fundamental question in Riemannian geometry. The classification in dimension $2$ is well-known and follows from the Gauss-Bonnet theorem. In dimension $3$, the classification   follows from Hamilton's work \cite{Ha}; in particular, a compact, simply connected  $3$-manifold of positive curvature must be diffeomorphic to the $3$-sphere. In dimension $n\geq 4$, in contrast,  a complete solution to the classification problem remains elusive to this day, as evidenced by the relative scarcity of examples and techniques for the construction of compact manifolds with positive or nonnegative curvature. 
Given these difficulties, it has been helpful to first consider the classification of the most symmetric spaces in these classes, that is, those  with a ``large'' group of isometries. This approach, proposed by Grove \cite{Gr2002}, allows for flexibility in deciding which isometry groups are to be considered ``large''. The classification of simply connected positively curved homogenenous spaces (cf.~\cite{Ber,Wa,AW,BeBe,Wi99}), for example,  may be framed in this program, which  has led to other classification results and to new examples of positively and nonnegatively curved Riemannian manifolds (cf.~\cite{GZ,GWZ08,Ve0, Ve,GVZ_GAFA,De}).

Let $(M,g)$ be a (compact) Riemannian manifold and and let $\Iso(M,g)$ be its isometry group, which is a  (compact) Lie group (cf.~\cite{MS,Ko}). There are several possible measures for the size of $\Iso(M,g)$, e.g.  the \emph{cohomogeneity}, defined as the dimension of the orbit space of the action of $\Iso(M,g)$ on $(M,g)$, the \emph{symmetry degree}, defined as the dimension of $\Iso(M,g)$, or the \emph{symmetry rank}, defined as the rank of $\Iso(M,g)$ and denoted by $\symrank(M,g)$. In this note we will focus our attention on this last invariant in the cases when $(M,g)$ has  positive curvature and when $(M,g)$ has \emph{quasipositive curvature},  i.e.~$(M,g)$ has nonngative curvature and a point with strictly positive curvature.

The following three problems arise naturally in the study of Riemannian manifolds and their symmetry rank:
\smallskip


\begin{itemize}
	\item[(a)] \emph{Maximal  symmetry rank:}  Given a class $\mathcal{M}^n$ of Riemannian $n$-manifolds, find an optimal upper bound $K$ for the symmetry rank of the elements in $\mathcal{M}^n$. \vspace{.1in}
	
	\item[(b)] \emph{Topological classification:}  Classify, up to diffeomorphism, all manifolds in $\mathcal{M}^n$ with symmetry rank $k\leq K$. \vspace{.1in}
	
	\item[(c)] \emph{Equivariant classification:} Let $M\in \mathcal{M}^n$ with $\symrank(M)=k$. Classify, up to equivariant diffeomorphism, all possible (effective) isometric actions of $T^k$ on $M$ and realize  these actions via appropriate Riemannian metrics on  $M$. 
\end{itemize} 
\smallskip

These problems have received particular attention when $\mathcal{M}^n$ is the class of compact, positively curved $n$-manifolds or the class of compact, simply connected $n$-manifolds of nonnegative curvature (cf.~\cite{GGS_2011_MSR,GGS_AMSRD5,GGK_T2D4T3D5,GroveSearle1994,HsKl1989,Kl1989,Ro2002,SeYa1994,Wi2003,Wi2006}). In a curvature-free setting,  analogs of problems (a), (b) and (c) for  compact, simply connected  smooth $n$-manifolds, $3\leq n\leq 6$, have also been extensively studied (cf.~\cite{F1,F2,Ko2006,N,Oh,Oh2,OrRa1967,OR,OR2,Ra1968}). 

The maximal symmetry rank problem and the topological classification of compact, positively curved manifolds of maximal symmetry rank were first considered by Grove and Searle \cite{GroveSearle1994}:


\begin{thm}[Grove, Searle \cite{GroveSearle1994}]
\label{T:MSR_PC} Let $(M^n,g)$ be a compact Riemannian $n$-manifold of positive curvature. Then the following hold:
	\begin{enumerate}
		\item $\operatorname{symrank}(M^n,g)\leq \lfloor (n+1)/2 \rfloor$.\\
		\item If $\operatorname{symrank}(M^n,g)=\lfloor (n+1)/2 \rfloor$, then $M^n$ is diffeomorphic to a sphere, a lens space or to a real or complex projective space. 
	\end{enumerate} 
\end{thm}

Let $(M,g_0)$ be isometric to any of the manifolds listed in Theorem~\ref{T:MSR_PC}-(2), equipped with its standard Riemannian metric $g_0$. As pointed out in \cite{GroveSearle1994},  $(M,g_0)$ has maximal symmetry rank. We will refer to the isometric torus actions on  $(M,g_0)$ as \emph{linear} torus actions.  Our first result is the equivariant classification of torus actions of maximal rank on compact, positively curved manifolds of maximal symmetry rank:


\begin{thm}
\label{T:EQC_PCMSR} Any effective, isometric torus action of maximal rank on a compact, positively curved Riemannian manifold of maximal symmetry rank is equivariantly diffeomorphic to a linear action.  
\end{thm}

It is natural to ask to what extent the conclusions of Theorem~\ref{T:MSR_PC} hold under weaker  curvature conditions, e.g. nonnegative curvature. In this case, an upper bound on the symmetry rank smaller than the dimension of the manifold, as in Theorem~\ref{T:MSR_PC},  cannot be achieved in full generality, since the $n$-dimensional flat torus has maximal symmetry rank $n$. Under the additional hypothesis of simple connectivity,  it has been conjectured (cf.~\cite{GGS_2011_MSR}) that if $(M^n,g)$ is a compact, simply connected nonnegatively curved Riemannian $n$-manifold, then $\symrank(M^n,g)\leq \lfloor 2n/3\rfloor$ and that, if  $n=3k$ and $(M^n,g)$ has maximal symmetry rank, then $M^n$ must be diffeomorphic to the product of $k$ copies of the $3$-sphere $\sphere^3$. This conjectural bound on the symmetry rank has been verified in dimensions at most $9$; the topological classification of compact, simply connected Riemannian manifolds of nonnegative curvature and maximal symmetry rank has also been completed in dimensions at most $6$, verifying the diffeomorphism conjecture in dimensions $3$ and $6$ (cf.~\cite{GGS_2011_MSR}).

In addition to nonnegatively curved (Riemannian) manifolds, one may consider manifolds with \emph{almost positive curvature}, i.e.~nonnegatively curved manifolds with positive curvature on an open and dense set, or manifolds with \emph{quasipositive curvature}, 
i.e.~nonnegatively curved manifolds with a point at which all tangent $2$-planes have positive curvature. These  two  families may be considered as intermediate classes  between positively and nonnegatively curved manifolds, and may be used as test cases to determine to what extent the collections of positively and nonnegatively curved manifolds  differ from each other. In the noncompact case, it follows from Perelman's proof of the Soul Conjecture \cite{Pe} that a complete, noncompact manifold with quasipositive curvature must be diffeomorphic to $\mathbb{R}^n$; in particular, it admits positive curvature. In the compact case, $\RP^2\times\RP^3$ admits a metric with quasipositive curvature (cf.~\cite{Wi02}) and cannot support a metric of positive curvature. In contrast to this, in the simply connected case there are no known obstructions distinguishing  compact manifolds with positive, almost positive, quasipositive or nonnegative curvature . 

Although there are many examples of manifolds with quasipositive or almost positive curvature (cf.~\cite{GM,PeWi,Fred01,Wi02,Ta,EK,Ke,DV, De2}),  including an exotic $7$-sphere,   the topological classification of these spaces remains open and one may consider problems (a), (b) and (c) for these  classes of Riemannian manifolds. 
Problem (a) was solved by Wilking \cite{Wi2008}, who showed that the bound for the symmetry rank  in  Theorem~\ref{T:MSR_PC}-(1) also holds for Riemannian manifolds with  quasipositive curvature:


\begin{thm}[Wilking \cite{Wi2008}] 
\label{T:MSR_QPC}
If $(M^n,g)$ is an $n$-dimensional Riemannian manifold of quasipositive curvature, then $\operatorname{symrank}(M^n,g)\leq \lfloor(n+1)/2\rfloor$.
\end{thm}

Our second result is the topological classification of  compact, simply-connected $4$- and $5$-manifolds of quasipositive curvature and maximal symmetry rank:%


\begin{thm}
\label{T:MSR_QPC_RIG}
Let $M^n$ be a compact, simply connected Riemannian  $n$-manifold with quasipositive curvature and maximal symmetry rank. Then the following hold:
	\begin{enumerate}
		\item If $n=4$, then $M^4$ is diffeomorphic to $\sphere^4$ or $\CP^2$.\\
		\item If  $n=5$, then $M^5$ is diffeomorphic to $\sphere^5$.   
	\end{enumerate}
\end{thm}
It follows from work of Orlik and Raymond \cite{OR}, in dimension $4$, and of Oh \cite{Oh}, in dimension $5$, that  smooth, effective $T^2$ actions on $\sphere^4$ and $\CP^2$, and  smooth, effective $T^3$ actions on $\sphere^5$, are equivalent to linear actions. Therefore, the isometric torus actions in Theorem~\ref{T:MSR_QPC_RIG} must be equivalent to linear actions; any such action can be realized via a standard metric that is  (trivially) quasipositively cuved, whence the equivariant classification follows.

Recall that the known examples of simply connected $4$- and $5$-manifolds of nonnegative curvature that are not known to admit positively curved metrics are $\sphere^2\times\sphere^2$, $\CP^2\#\pm\CP^2$, $\sphere^2\times \sphere^3$, the non-trivial bundle $\sphere^2\tilde{\times}\sphere^3$ and the Wu manifold $\SU(3)/\SO(3)$. 
 Out of these manifolds, only $\sphere^2\times\sphere^2$ and $\CP^2\#\pm\CP^2$ admit metrics with nonnegative curvature and maximal symmetry rank $2$, and the bundles $\sphere^2\times \sphere^3$ and $\sphere^2\tilde{\times}\sphere^3$ are the only ones  admitting metrics of nonnegative curvature and maximal symmetry rank $3$; the Wu manifold $\SU(3)/\SO(3)$, equipped with its standard nonnegatively curved homgoeneous metric, has symmetry rank $2$ (cf.~\cite{GGS_AMSRD5,GGK_T2D4T3D5}). On the other hand, the trivial sphere bundle $\sphere^3\times\sphere^2$ carries an almost positively curved metric with symmetry rank $1$ (cf.~\cite{Wi02}) and it is not known if the remaining $4$- and $5$-manifolds  listed in this paragraph admit metrics of quasipositive curvature. Theorem~\ref{T:MSR_QPC_RIG} implies that any such metric would have symmetry rank at most $1$, in dimension $4$, and at most $2$, in dimension $5$.

We conclude these remarks by recalling the so-called \emph{deformation conjecture} (cf.~\cite{GM,Wi02}), which states that if $(M,g)$ is a complete Riemannian manifold of quasipositive curvature, then $M$ admits a metric with positive curvature. As pointed out above, this conjecture  is true if $M$ is noncompact, false  if $M$ is compact and not simply connected, and remains open if $M$ is compact and simply connected (see \cite{PW} for the construction of a metric with positive curvature on the Gromoll-Meyer sphere, an exotic $7$-sphere with quasipositive curvature). Theorem~\ref{T:MSR_QPC_RIG} may be seen as supporting this conjecture when $(M,g)$ is compact, simply connected and has maximal symmetry rank.


The contents of this note are organized as follows. In section~\ref{S:Preliminaries} we  collect some background material and recall the proof of Theorem~\ref{T:MSR_QPC}. This result was not available in the literature; for the sake of reference, we have included Wilking's proof as conveyed to us by M. Kerin.  
In section~\ref{S:pc_msr_equiv_classif} we prove Theorem~\ref{T:EQC_PCMSR}  and in section~\ref{S:QPC_symrk_topclas_d4} we prove Theorem ~\ref{T:MSR_QPC_RIG}. The proofs follow easily from restrictions on the structure of the manifolds and their orbit spaces imposed by the curvature hypotheses and the rank of the actions.  As the reader may have already noticed, we have strived to give extensive references to the literature. 

\begin{ack} This note originated from talks I gave at the \emph{Tercer Miniencuentro de Geometr\'ia Diferencial} in December of 2010,  at CIMAT, M\'exico. It is a pleasure to thank Rafael Herrera, Luis Hern\'andez Lamoneda and CIMAT for their hospitality during the workshop. I also  wish to thank Martin Kerin, for  sharing his notes with the proof of Theorem~\ref{T:MSR_QPC}, and Xiaoyang Chen, for some conversations on quasipositively curved manifolds. 
\end{ack}


\section{Preliminaries}
\label{S:Preliminaries}


\subsection{Basic setup and notation.}
Let $G\times M \rightarrow M$,  $\ m \mapsto g(m)$, be a smooth action of a compact Lie group $G$ on a smooth manifold $M$. The orbit $G(p)$ through a point $p \in M$ is diffeomorphic to the quotient $G/ G_p$, where $G_p = \{g \in G\, :\, g(p) = p \}$ is the \emph{isotropy} subgroup of $G$ at $p$.  If $G_p$ acts trivially on the normal space to the orbit at $p$, then $G / G_p$ is called a \emph{principal orbit}.  The set of principal orbits is open and dense in $M$. Since the isotropy groups of principal orbits are all conjugate in $G$, all principal orbits have the same dimension. The isotropy group of principal orbits is the \emph{principal isotropy subgroup}. If $G/G_p$ has the same dimension as a principal orbit and $G_p$ acts non-trivially on the normal space at $p$, then $G / G_p$ is called an \emph{exceptional orbit}.  An orbit that is neither principal nor exceptional is called a \emph{singular orbit}.  When $G_p = G$, the point $p$ is called a \emph{fixed point} of the action.   Recall that the \emph{ineffective kernel} of the action  is $K := \{g \in G \, : \ g(m) = m, \ \text{for all} \ m \in M\}$.  The action is \emph{effective} if the ineffective kernel is trivial.  The group $\widetilde G = G/K$  always acts effectively on $M$.

Given a subset $X \subset M$, we will denote its projection under the orbit map $\pi: M \rightarrow M/G$ by $X^*$. Following this convention, we will denote the orbit space $M/G$  by $M^*$. 

Recall that a finite dimensional length space $(X,\mathrm{dist})$ is an \emph{Alexandrov space} if it has curvature bounded from below in the triangle comparison sense (cf.~\cite{BBI}). When $(M,g)$ is a complete, connected Riemannian manifold and $G$ is a compact Lie group acting on $(M,g)$ by isometries, the orbit space $M^*$ can be made into a metric space $(M^*,\mathrm{dist})$  by defining the distance between orbits $p^*$ and $q^*$ in $M^*$ as the distance between the orbits $G(p)$ and $G(q)$ as subsets of $(M,g)$. If, in addition, $(M,g)$ has sectional curvature bounded below by $k$, then the orbit space $(M^*,\mathrm{dist})$ equipped with this so-called \emph{orbital metric} is an Alexandrov space with curvature bounded below by $k$. The \emph{space of directions} of a general Alexandrov space at a point $x$ is,
by definition, the completion of the 
space of geodesic directions at $x$. In the case of an orbit space $M^*=M/G$, the space of directions $\Sigma_{p^*}M^*$ at a point $p^*\in M^*$ consists of geodesic directions and is isometric to
\[
\sphere^{\perp}_p/G_p,
\] 
where $\sphere^{\perp}_p$ is the unit normal sphere to the orbit $G(p)$ at $p\in M$.


\subsection{Proof of Theorem~\ref{T:MSR_QPC} (Wilking \cite{Wi2008})}
\label{SS:MSR_QPC}

Let $(M^n,g)$ be an $n$-dimensional Riemannian manifold of quasipositive curvature with an (effective) isometric  $T^k$ action.  It suffices to show that if $k>(n+1)/2$, then $(M^n,g)$ cannot have quasipositive curvature. Throughout the proof we will let $\Omega_p = T^k(p)$ be a principal orbit of the $T^k$ action for some  $p\in M$. Given $q\in \Omega_p$, we  let $T_q(\Omega_p)^\perp$ be the orthogonal complement of $T_q(\Omega_p)$ in $T_pM$. Recall that the second fundamental form at $q\in \Omega_p$ is given by 
\[
\alpha:T_q(\Omega_p)\times T_q(\Omega_p)\rightarrow T_q(\Omega_p)^\perp.
\] 

Let $u\in T_q(\Omega_p)$, $|u|=1$, and let $u^\perp$ be its orthogonal complement in $T_q(\Omega_p)$. Note that
\[
\dim u^{\perp} = \dim \Omega_p -1 > \frac{n-1}{2},
\]
\[
\dim T_q(\Omega_p)^\perp = \dim M - \dim \Omega_p < \frac{n-1}{2}.
\]
Consider $\alpha(u,\cdot):u^\perp \rightarrow  T_q(\Omega_p)^\perp$. For dimension reasons there exists a unit vector $w\in u^\perp$ such that $\alpha(u,w)=0$.

We will now show that there exists a unit vector $v\in T_q(\Omega_p)$ such that $\alpha(v,v)=0$. Suppose that there is no such vector in $T_q\Omega_p$. Then there exists $u\in T_q(\Omega_p)$, $|u|=1$, such that $|\alpha(u,u)|>0$ is minimal. By the preceding paragraph, there exists $w\in u^\perp$, $|w|=1$, such that $\alpha(u,w)=0$. Consider the function 
\begin{eqnarray*}
f(t) 	& :=  & |\alpha((\cos t) u + (\sin t) w, (\cos t) u + (\sin t) w)|^2\\
	& =  & |(\cos^2 t)\alpha(u,u)+(\sin^2 t)\alpha(w,w)|^2.
\end{eqnarray*}

Since $f(0)=|\alpha(u,u)|^2$ is minimal, $f'(0)=0$ and 
\[
0\leq f''(0)=4(\langle\,\alpha(u,u),\alpha(w,w)  \,\rangle - |\alpha(u,u)|^2).
\]
In particular, since $|\alpha(u,u)|^2>0$, we have that
\[
\langle\, \alpha(u,u),\alpha(w,w)\, \rangle >0.
\]
It then follows from the Gauss formula that
\[
\sec_{\Omega_p}(u,w)=\sec_M(u,w)+\langle\, \alpha(u,u),\alpha(w,w)\,\rangle >0.
\]
This yields a contradiction, since $\Omega_p$ is a torus equipped with a left-invariant metric, hence flat. 

It follows from the preceding paragraph that there exist orthogonal unit vectors $u,v\in T_q(\Omega_p)$ such that $\alpha(u,u)=\alpha(u,v)=0$. Then, by the Gauss formula, 
\[
\sec_{\Omega_p}(u,v)=\sec_M(u,v).
\]
Since the choice of principal orbit $\Omega_p$ and $q\in \Omega_p$ was arbitrary, it follows that there is an open and dense set of points $p\in M$ with a tangent plane $\Pi_p\subset T_p(\Omega_p)$ such that 
\[
0=\sec_{\Omega_p}(\Pi_p) = \sec_M (\Pi_p).
\]
Therefore, $(M,g)$ cannot have quasipositive curvature. \hfill $\square$


\section{Proof of Theorem~\ref{T:EQC_PCMSR}}
\label{S:pc_msr_equiv_classif}

Let $(M^n,g)$ be a compact, positively curved Riemannian $n$-manifold of maximal symmetry rank with an (effective) isometric action of a torus $T^k$ of maximal rank.  We first recall some basic facts from \cite{GroveSearle1994}. There exists  a circle subgroup $T^1\leq T^k$ with fixed point set a totally geodesic codimension $2$ submanifold $F^{n-2}$ of $M^n$ such that $F^{n-2}=\partial (M^n/T^1)$. The orbit space $M^n/T^1$ is a positively curved Alexandrov space homeomorphic to the cone over an orbit $p^*$ at maximal distance from $\partial(M^n/T^1)$. The  isotropy subgroup of $p^*$ is either $T^1$, $\mathbb{Z}_k$, $k\geq 2$, or $\mathsf{1}$, the trivial subgroup of $T^1$. Since $F^{n-2}$ is diffeomorphic to the space of directions of $M^n/T^1$ at $p^*$,  $F^{n-2}$ is diffeomorphic to a sphere, if $p^*$ is a principal orbit; a lens space or an even-dimensional real projective space, if $p^*$ is an exceptional orbit; or to a complex projective space, if $p^*$ is a fixed point. Moreover, there exists an invariant disc bundle decomposition 
\[
M^n=D(F^{n-2})\cup_E D(G(p)),
\]
where $D(F^{n-2})$ is a  tubular neighborhood of $F^{n-2}$, $D(G(p))$ is a tubular neighborhood of the orbit $G(p)$ corresponding to the vertex $p^*$ of the orbit space, and $E$ is the common boundary $\partial D(F^{n-2})=\partial D(G(p))$. The manifold $M^n$ is diffeomorphic to
 a sphere if $p^*$ is a principal orbit;  a lens space or a real projective space, if $p^*$ is an exceptional orbit;  or to  a complex projective space, if $p^*$ is a fixed point.

We will now prove the theorem in the case where $M^n$ is diffeomorphic to an $n$-sphere $\sphere^n$. We proceed by induction on the dimension $n$. For $n= 2$, it is well known that any smooth $T^1$ action on $\sphere^2$ is equivalent to a linear action. Fix $n>2$ and let $(\sphere^n,g)$ be a positively curved $n$-sphere of maximal symmetry rank, so that there exists an effective isometric  $T^k$ action on $(\sphere^n,g)$ with $k= \lfloor (n+1)/2 \rfloor$. As recalled in the preceding paragraph, there is a circle subgroup $T^1\leq T^k$ with fixed point set a totally geodesic sphere  $\sphere^{n-2}\subset \sphere^n$ of  codimension $2$. The invariant decomposition of $\sphere^n$ into a union of disc bundles induced by the $T^1$ action is given by
\begin{eqnarray*}
	\sphere^n		& 	\simeq	&	D(\sphere^{n-2})\cup_{\partial D(\sphere^1)} D(\sphere^1)\\[.1in]
				&	\simeq	& 	(\sphere^{n-2}\times D^2) \cup_{\sphere^{n-2}\times \sphere^1} (D^{n-1}\times \sphere^1),
\end{eqnarray*}
where $D(\sphere^{n-2})$ is a  tubular neighborhood of the fixed point set $\sphere^{n-2}$ of $T^1$ and  $D(\sphere^1)$ is a tubular neighborhood of the orbit $T^1(p)\simeq \sphere^1$ whose projection $p^*$ is the vertex point of the orbit space $\sphere^n/T^1$. As in \cite{GroveSearle1994}, the $T^1$ action on 
\[
D(\sphere^{n-2}) \simeq \sphere^{n-2}\times D^2
\]
is equivalent to the $T^1$ action on  $ \sphere^{n-2}\times_{T^1} D^2$, the associated disc bundle to the $T^1$ action on the (trivial) normal bundle of the fixed point set $\sphere^{n-2}$.

We may write $T^k=T^1\oplus T^{k-1}$, where $T^{k-1}$ is the orthogonal complement of $T^1$ in $T^k$. Observe now that $T^{k-1}$ acts effectively and isometrically on $\sphere^{n-2}\subset \sphere^n$. By induction, the action of $T^{k-1}$ on $\sphere^{n-2}$ is linear. It follows that the $T^k$ action on $\sphere^{n-2}\times D^2$ is given by the product of a linear $T^{k-1}$ action on $\sphere^{n-2}$ and a linear $T^1$ action on $D^2$. Consequently, on the boundary $ \sphere^{n-2}\times \sphere^1$ the $T^k$ action is the product of a linear $T^{k-1}$ action on $\sphere^{n-2}$ and a linear $T^1$ action on $\sphere^1$.  On $D^{n-1}\times \sphere^1$, the other half of the disc bundle decomposition of $\sphere^n$, the $T^k$ action is given by the product of a linear $T^{k-1}$ action on $D^{n-1}$ and a linear action $T^1$ action on $\sphere^1$. Observe that the linear $T^{k-1}$ action on $D^{n-1}$ is the cone over the linear action of $T^{k-1}$ on the $\sphere^{n-2}$ factor of the boundary $D(\sphere^1)\simeq \sphere^{n-2}\times \sphere^1$. Hence, the $T^k$ action on  $(\sphere^n,g)$ is equivariantly diffeomorphic to the linear $T^k$ action on $\sphere^n=\sphere^{1}*\sphere^{n-2}\subset \mathbb{R}^{2}\times\mathbb{R}^{n-2}$ given by letting $T^1$ act orthogonally on  $\mathbb{R}^{2}$ and $T^{k-1}$ act orthogonally on $\mathbb{R}^{n-2}$.

When $M^n$ is diffeomorphic to  a lens space or to a real projective space, the conclusion follows by passing to the universal covering space and observing that  the covering torus action must be equivalent to a linear action on $\sphere^n$.   

The proof when $M^n$ is diffeomorphic to  $\CP^m$ is analogous to to the case of the sphere. For $m\geq 2$, the equivariant disc bundle decomposition is given by 
\[
\CP^{m} \simeq \sphere^{2n-1}\times_{T^1}D^2\cup_{\sphere^{2n-1}} D^{2n},
\]
 where $T^1$ is the circle subgroup of $T^m$ fixing both $\CP^{m-1}\subset \CP^{m}$ and the vertex of $D^{2n}$, and  $\sphere^{2m-1}\times_{T^1}D^2$ is the normal disc bundle of $\CP^{m-1}$ in $\CP^{m}$. The $T^m$ action is equivalent to a linear  $T^m$ action on $\CP^m$ induced by  a linear $T^{m+1}$ action on $\sphere^{2m+1}$ via the projection map $\pi: \sphere^{2m+1}\rightarrow \CP^m$ of the Hopf action. \hfill $\square$


\section{Proof of Theorem~\ref{T:MSR_QPC_RIG}}
\label{S:QPC_symrk_topclas_d4}
We proceed along the lines of  \cite{GGS_2011_MSR}. Let $(M^n,g)$ be a compact, simply connected Riemannian $n$-manifold, $n=4$ or $5$, with quasipositive curvature and maximal symmetry rank.
Then $(M^n,g)$ has an isometric torus action whose orbit space $M^*$ is $2$-dimensional. It follows from work of several authors (cf.~\cite{Br1972,KMP,Oh,OR}) that the orbit space $M^*$ of the action has the following properties: $M^*$ is homeomorphic to   a $2$-disk, the boundary of $M^*$ is the set of singular orbits and the interior of $M^*$ consists of principal orbits. Moreover, when $n=4$, there are at least two isolated orbits with isotropy $T^2$ and, when $n=5$, there are at least three isolated orbits with isotropy $T^2$. In both cases, points in the the arcs  in the boundary of $M^*$ joining  orbits with isotropy $T^2$ have isotropy conjugate to a circle $T^1$; the angle between these arcs is $\pi/2$ and is the length of the space of directions at an orbit with isotropy $T^2$ in $M^*$

Since $(M^n,g)$ is a quasipositively curved Riemannian manifold,  $M^*$ is a nonnegatively curved $2$-manifold with non-smooth boundary and positive curvature on an open subset. A simple comparison argument using  Toponogov's theorem  shows that there can be at most $4$ points in $M^*$ corresponding to orbits with isotropy $T^2$ and, if there are  $4$ such points, then $M^*$ must be isometric to a flat rectangle. Since $M^*$ has positive curvature on an open subset,  there can be at most $3$ points in $M^*$ with isotropy $T^2$. Hence, for $n=4$, the orbit space $M^*$ has $2$ or $3$ points with isotropy $T^2$ and, for $n=5$, $M^*$ has exactly $3$ such points. The conclusions of the theorem now follow from the Orlik-Raymond classification of compact, smooth, simply connected $4$-manifolds with a smooth, effective $T^2$ action (cf.~\cite{OR}), and from Oh's classification of compact, smooth, simply connnected $5$-manifolds with a smooth, effective $T^3$ action (cf.~\cite{Oh}). 
\hfill $\square$

\bibliographystyle{amsplain}


\end{document}